\documentclass[runningheads]{asmda2005}
\usepackage{amssymb}

\usepackage{asmda2005References} 
\newtheorem{Lem}{\underline {Lemma}}

\newtheorem{Theorem}{\underline {Theorem}}

\begin{document}
%
\title*{Testing the number of parameters of multidimensional MLP}
%
\toctitle{Testing the number of parameters with multidimensional MLP}
%
\titlerunning{multidimensional MLP}
%
\author{
  Joseph Rynkiewicz\inst{1}
}
%
\index{Rynkiewicz, J.}
%
\authorrunning{Joseph Rynkiewicz}
%
\institute{
  SAMOS - MATISSE \\
  Universit\'e de Paris I\\
  72 rue Regnault, 75013 Paris, France\\
  (e-mail: {\tt joseph.rynkiewicz@univ-paris1.fr})
}

\maketitle             

\begin{abstract}
This work concerns testing the number of parameters in one hidden layer multilayer perceptron (MLP). For this purpose we assume that we have identifiable models, up to a finite group of transformations on the  weights, this is for example the case when the number of hidden units is know. In this framework, we show that we get a simple asymptotic distribution, if we use the logarithm of the determinant of the empirical error covariance matrix as cost function. 
\end{abstract}

\section{Introduction}
Consider a sequence $\left(Y_{t},Z_{t}\right)_{t\in \mathbb{N}}$
of i.i.d.\footnote{It is not hard to extend all what we show in this paper for stationary mixing variables and so for time series} random vectors (i.e. identically distributed and independents). So, each couple $\left(Y_{t},Z_{t}\right)$ has the same law that a generic variable $(Y,Z)\in\mathbb R^d\times\mathbb R^{d'}$. 
\subsection{The model}
Assume that the model can be written 
\[
Y_{t}=F_{W^0}(Z_{t})+\varepsilon _{t}
\]
where 

\begin{itemize}
\item $F_{W^0}$ is a function represented by a one hidden layer MLP with parameters or
weights $W^0$ and sigmoidal functions in the hidden unit. 
\item The noise, $(\varepsilon _{t})_{t\in \mathbb N}$, is sequence of i.i.d. centered variables with unknown invertible covariance matrix $\Gamma(W^0)$. Write $\varepsilon$ the generic variable with the same law that each $\varepsilon_t$.
\end{itemize}
Notes that a finite number of transformations of the weights leave the MLP functions invariant, these permutations form a finite group (see \cite{Sussman}). To overcome this problem, we will consider equivalence classes of MLP : two MLP are in the same class if the first one is the image by such transformation of the second one, the considered set of parameter is then the quotient space of parameters by the finite group of transformations. 

In this space, we assume that the model is identifiable, this can be done if we consider only MLP with the true number of hidden units (see \cite{Sussman}). Note that, if the number of hidden units is over-estimated, then such test can have very bad behavior (see \cite{Fukumizu}). We agree that the assumption of identifiability is very restrictive, but we want emphasize the fact that, even in this framework, classical test of the number of parameters in the case of multidimensional output MLP is not satisfactory and we propose to improve it. 
\subsection{testing the number of parameters}
Let $q$ be an integer lesser than $s$, we want to test ``$H_0 : W\in \Theta_q \subset \mathbb R^q$'' against  ``$H_1 : W\in \Theta_s \subset \mathbb R^s$'', where the sets $\Theta_q$ and $\Theta_s$ are compact. $H_0$ express the fact that $W$ belongs to a subset of $\Theta_s$ with a parametric dimension lesser than $s$ or, equivalently, that $s-q$ weights of the MLP in $\Theta_s$ are null. 
If we consider the classic cost function : $V_n(W)=\sum_{t=1}^n\Vert Y_{t}-F_{W}(Z_{t})\Vert^2$ where $\Vert x\Vert$ denotes the Euclidean norm of $x$, we get the following statistic of test : 
\[
S_n=n\times\left(\min_{W\in \Theta_q}V_n(W)-\min_{W\in \Theta_s}V_n(W)\right)
\]
It is shown in \cite{Yao}, that $S_n$ converges in law to a ponderated sum of $\chi^2_1$
\[
S_n\stackrel{\cal D}{\rightarrow}\sum_{i=1}^{s-q}\lambda_i\chi_{i,1}^2
\]
where the $\chi_{i,1}^2$ are $s-q$ i.i.d. $\chi^2_1$ variables and $\lambda_i$ are strictly positives values, different of 1 if the true covariance matrix of the noise is not the identity. So, in the general case, where the true covariance matrix of the noise is not the identity, the asymptotic distribution is not known, because the $\lambda_i$ are not known and it is difficult to compute the asymptotic level of the test. 

To overcome this difficulty we propose to use instead the cost function
\begin{equation}
U_{n}\left(W\right):=\ln \det \left(\frac{1}{n}\sum _{t=1}^{n}(Y_{t}-F_{W}(Z_{t}))(Y_{t}-F_{W}(Z_{t}))^{T}\right).
\end{equation}
we will show that, under suitable assumptions, the statistic of test :
\begin{equation}\label{tn}
T_n=n\times\left(\min_{W\in \Theta_q}U_n(W)-\min_{W\in \Theta_s}U_n(W)\right)
\end{equation}
will converge to a classical $\chi^2_{s-q}$ so the asymptotic level of the test will be very easy to compute. The sequel of this paper is devoted to the proof of this property.

\section{Asymptotic properties of $T_n$}
In order to investigate the asymptotic properties of the test we have to prove the consistency and the asymptotic normality of $\hat W_n=\arg\min_{W\in\Theta_s}U_n(W)$.
Assume, in the sequel, that $\varepsilon$ has a moment of order at least 2 and note 
\[
\Gamma_n(W)=\frac{1}{n}\sum _{t=1}^{n}(Y_{t}-F_{W}(Z_{t}))(Y_{t}-F_{W}(Z_{t}))^{T}
\]  
remark that these matrix $\Gamma_{n}(W)$ and it inverse  are symmetric.
in the same way, we note $\Gamma(W)=\lim_{n\rightarrow\infty}\Gamma_n(W)$, which is well defined because of the moment condition on $\varepsilon$
\subsection{Consistency of $\hat W_n$}
First we have to identify contrast function associated to $U_n(W)$ 
\begin{Lem}
\[
U_n(W)-U_n(W^0)\stackrel{a.s.}{\rightarrow}K(W,W^0)
\]
with $K(W,W^0)\geq 0$ and $K(W,W^0)= 0$ if and only if $W=W^0$.
\end{Lem}
\paragraph{Proof : }
By the strong law of large number we have
\[
\begin{array}{l}
U_n(W)-U_n(W^0)\stackrel{a.s.}{\rightarrow}\ln \det(\Gamma(W))-\ln \det(\Gamma(W^0))=\ln\frac{\det(\Gamma(W))}{\det(\Gamma(W^0))}=\\
\ln \det\left(\Gamma^{-1}(W^0) \left(\Gamma(W)-\Gamma(W^0)\right)+I_d\right)
\end{array}
\]
where $I_d$ denotes the identity matrix of $\mathbb R^d$.
So,  the lemme is true if $\Gamma(W)-\Gamma(W^0)$ is a positive matrix, null only if $W=W^0$. But this property is true since
\[
\begin{array}{l}
\Gamma(W)=E\left((Y-F_W(Z))(Y-F_W(Z))^T\right)=\\
E\left((Y-F_{W^0}(Z)+F_{W^0}(Z)-F_W(Z))(Y-F_{W^0}(Z)+F_{W^0}(Z)-F_W(Z))^T\right)=\\
E\left((Y-F_{W^0}(Z))(Y-F_{W^0}(Z))^T\right)+\\
E\left((F_{W^0}(Z)-F_W(Z))(F_{W^0}(Z)-F_W(Z))^T\right)=\\
\Gamma(W^0)+E\left((F_{W^0}(Z)-F_W(Z))(F_{W^0}(Z)-F_W(Z))^T\right)\blacksquare
\end{array}
\]
We deduce then the theorem of consistency :
\begin{Theorem}
If $E\left( \Vert \varepsilon \Vert^2\right)<\infty$, 
\[
\hat W_n\stackrel{P}{\rightarrow}W^0
\]
\end{Theorem} 
\paragraph{Proof} Remark that it exist a constant $B$ such that 
\[
sup_{W\in\Theta_s}\Vert Y-F_W(Z)\vert^2<\Vert Y\Vert^2+B
\]
because $\Theta_s$ is compact, so $F_W(Z)$ is bounded. 
For a matrix $A\in\mathbb R^{d\times d}$, let $\Vert A\Vert$ be a norm, for example $\Vert A\Vert^2=tr\left(AA^T\right)$. We have
\[
\begin{array}{l}
\lim \inf_{W\in \Theta_s} \Vert \Gamma_n(W)\Vert =\Vert \Gamma(W^0)\Vert>0\\
\lim \sup_{W\in \Theta_s} \Vert \Gamma_n(W)\Vert :=C <\infty
\end{array}
\]
and since the function :
\[
\Gamma\mapsto \ln \det \Gamma, \mbox{ for }C\geq\Vert\Gamma\Vert\geq\Vert\Gamma(W^0)\Vert
\] is uniformly continuous, by the same argument that example 19.8 of \\
 \cite{Vandervaart} the set of functions $U_n(W),\ W\in\Theta_s$ is Glivenko-Cantelli.

Finally, the theorem 5.7 of \cite{Vandervaart}, show that $\hat W_n$ converge in probability to $W^0$ $\blacksquare$.
\subsection{Asymptotic normality}
 For this purpose we have to compute the first and the second derivative with respect to the parameters of $U_n(W)$. First, we introduce a notation :  if $F_W(X)$ is a $d$-dimensional parametric function depending of a parameter $W$, write $\frac{\partial F_W(X)}{\partial W_k}$ (resp. $\frac{\partial^2 F_W(X)}{\partial W_k\partial W_l}$) for the $d$-dimensional vector of partial derivative (resp. second order partial derivatives) of each component of $F_W(X)$. 
\paragraph{First derivatives : }
 if $\Gamma_n(W)$ is a matrix depending of the parameter vector $W$, we get  from  \cite{Magnus} \[
\frac{\partial }{\partial W_{k}}\ln \det \left(\Gamma _{n}(W)\right)=tr\left(\Gamma_{n}^{-1}(W) \frac{\partial }{\partial W_{k}}\Gamma _{n}(W)\right)\]

Hence, if we note 
\[
A_n(W_k)=\frac{1}{n}\sum _{t=1}^{n}\left(-\frac{\partial F_{W}(z_{t})}{\partial W_{k}}(y_t-F_{W}(z_{t}))^T\right)
\]
using the fact 
\[
tr\left(\Gamma_{n}^{-1}(W)A_n(W_k)\right)=tr\left(A_n^T(W_k)\Gamma_{n}^{-1}(W)\right)=tr\left(\Gamma_{n}^{-1}(W)A_n^T(W_k)\right)
\]

we get 
\begin{equation}\label{first_deriv}
\frac{\partial }{\partial W_{k}}\ln \det \left(\Gamma _{n}(W)\right)=2tr\left(\Gamma_{n}^{-1}(W)A_n(W_k)\right)
\end{equation}
\paragraph{Second derivatives : }
We write now 
\[
B_n(W_k,W_l):=\frac{1}{n}\sum _{t=1}^{n}\left( \frac{\partial F_{W}(z_{t})}{\partial W_{k}}\frac{\partial F_{W}(z_{t})}{\partial W_{l}}^T\right)
\]
and
\[
C_n(W_k,W_l):=\frac{1}{n}\sum _{t=1}^{n}\left( -(y_t-F_{W}(z_{t})) \frac{\partial^2 F_{W}(z_{t})}{\partial W_{k}\partial W_{l}}^T\right)
\]
We get

\[
\begin{array}{l}
\frac{\partial^2 U_n(W)}{\partial W_k\partial W_l}=\frac{\partial }{\partial W_{l}}2tr\left(\Gamma_{n}^{-1}(W)A_n(W_k)\right)=\\
2tr\left(\frac{\partial\Gamma_{n}^{-1}(W)}{\partial W_{l}}A_n(W_k)\right)+2tr\left(\Gamma_{n}^{-1}(W)B_n(W_k,W_l)\right)+2tr\left(\Gamma_{n}(W)^{-1}C_n(W_k,W_l)\right)\\
\end{array}
\]
Now, \cite{Magnus}, give an analytic form of the derivative of an inverse matrix, so we get
\[
\begin{array}{l}
\frac{\partial^2 U_n(W)}{\partial W_k\partial W_l}=2tr\left(\Gamma_{n}^{-1}(W)\left(A_n(W_k)+A_n^T(W_k)\right)\Gamma_{n}^{-1}(W)A_n(W_k)\right)+\\
2tr\left(\Gamma_{n}^{-1}(W)B_n(W_k,W_l)\right)+2tr\left(\Gamma_{n}^{-1}(W)C_n(W_k,W_l)\right)
\end{array}
\]
so
\begin{equation}\label{second_deriv}
\begin{array}{l}
\frac{\partial^2 U_n(W)}{\partial W_k\partial W_l}=4tr\left(\Gamma_{n}^{-1}(W)A_n(W_k)\Gamma_{n}^{-1}(W)A_n(W_k)\right)\\
+2tr\left(\Gamma_{n}^{-1}(W)B_n(W_k,W_l)\right)+2tr\left(\Gamma_{n}^{-1}(W)C_n(W_k,W_l)\right)
\end{array}
\end{equation}
\paragraph{Asymptotic distribution of $\hat W_n$ :}
The previous equations allow us to give the asymptotic properties of the estimator minimizing the cost function $U_n(W)$, namely from equation (\ref{first_deriv}) and  (\ref{second_deriv}) we can compute the asymptotic properties of the first and the second derivatives of $U_n(W)$. If the variable $Z$ has a moment of order at least 3 then  we get the following lemma :
\begin{Theorem}

Assume that $E\left(\Vert\varepsilon\Vert^2\right)<\infty$ and $E\left(\Vert Z\Vert^3\right)<\infty$, let $\Delta U_n(W^0)$ be the gradient vector of $U_n(W)$ at $W^0$ and $HU_n(W^0)$ be the Hessian matrix of  $U_n(W)$ at $W^0$. 

Write finally
\[
B(W_k,W_l):=\frac{\partial F_{W}(Z)}{\partial W_{k}}\frac{\partial F_{W}(Z)}{\partial W_{l}}^T
\]
We get then 
\begin{enumerate}
\item $HU_n(W^0)\stackrel{a.s.}{\rightarrow}2I_0$
\item $\sqrt{n}\Delta U_n(W^0)\stackrel{Law}{\rightarrow}{\cal N}(0,4I_0)$
\item $\sqrt{n}\left(\hat W_n-W^0\right)\stackrel{Law}{\rightarrow}{\cal N}(0,I^{-1}_0)$
\end{enumerate}
where, the component $(k,l)$ of the matrix $I_0$ is :
\[
tr\left(\Gamma^{-1}_0E\left(B(W^0_k,W^0_l)\right) \right)
\]
\end{Theorem}
\paragraph{proof : }
We can show easily that, for all $x\in\mathbb R^d$, we have :
\[
\begin{array}{l}
\Vert \frac{\partial F_W(Z)}{\partial W_k}\Vert\leq Cte(1+\Vert Z\Vert)\\
\Vert \frac{\partial^2 F_W(Z)}{\partial W_k\partial W_l}\Vert\leq Cte(1+\Vert Z\Vert^2)\\
\Vert \frac{\partial^2 F_W(Z)}{\partial W_k\partial W_l}-\frac{\partial^2 F_W^0(Z)}{\partial W_k\partial W_l}\Vert\leq Cte\Vert W-W^0 \Vert(1+\Vert Z\Vert^3)
\end{array}
\]
Write
\[
A(W_k)=\left(-\frac{\partial F_{W}(Z)}{\partial W_{k}}(Y-F_{W}(Z))^T\right)
\]
and $U(W):=\log\det(Y-F_W(Z))$.

Note that the component $(k,l)$ of the matrix $4I_0$ is:
\[
E\left(\frac{\partial U(W^0)}{\partial W_k}\frac{\partial U(W^0)}{\partial W^0_l}\right)=E\left(2tr\left(\Gamma^{-1}_0A^T(W^0_k)\right)\times2tr\left(\Gamma^{-1}_0A(W^0_l)\right)\right)
\]
and, since the trace of the product is invariant by circular permutation, 
\[
\begin{array}{l}
E\left(\frac{\partial U(W^0)}{\partial W_k}\frac{\partial U(W^0)}{\partial W^0_l}\right)=\\
4E\left( -\frac{\partial F_{W^0}(Z)^T}{\partial W_k}\Gamma^{-1}_0(Y-F_{W^0}(Z))(Y-F_{W^0}(Z))^T\Gamma^{-1}_0\left(-\frac{\partial F_{W^0}(Z))}{\partial W_l}\right)\right)\\
=4E\left(\frac{\partial F_{W^0}(Z)^T}{\partial W_k}\Gamma^{-1}_0\frac{\partial F_{W^0}(Z)}{\partial W_l}\right)\\
=4tr\left(\Gamma^{-1}_0E\left(\frac{\partial F_{W^0}(Z)}{\partial W_k}\frac{\partial F_{W^0}(Z)^T}{\partial W_l}\right) \right)\\
=4tr\left(\Gamma^{-1}_0E\left(B(W^0_k,W^0_l)\right)\right) 

\end{array}
\]
Now, the derivative $\frac{\partial F_W(Z)}{\partial W_k}$ is square integrable, so $\Delta U_n(W^0)$ fulfills Lindeberg's condition (see \cite{Hall}) and 
\[
\sqrt{n}\Delta U_n(W^0)\stackrel{Law}{\rightarrow}{\cal N}(0,4I_0)
\] 
For the component $(k,l)$ of the expectation of the Hessian matrix, remark first that 
\[
\lim_{n\rightarrow \infty}tr\left(\Gamma_{n}^{-1}(W^0)A_n(W^0_k)\Gamma_{n}^{-1}(W^0)A_n(W^0_k)\right)=0
\]
and
\[
\lim_{n\rightarrow \infty}tr\Gamma_{n}^{-1}C_n(W^0_k,W^0_l)=0
\]
so
\[
\begin{array}{l}
\lim_{n\rightarrow \infty}H_n(W^0)=\lim_{n\rightarrow \infty}4tr\left(\Gamma_{n}^{-1}(W^0)A_n(W^0_k)\Gamma_{n}^{-1}(W^0)A_n(W^0_k)\right)+\\
2tr\Gamma_{n}^{-1}(W^0)B_n(W^0_k,W^0_l)+2tr\Gamma_{n}^{-1}C_n(W^0_k,W^0_l)=\\
=2tr\left(\Gamma^{-1}_0E\left(B(W^0_k,W^0_l)\right)\right)\\
\end{array}
\]
Now, since $\Vert \frac{\partial^2 F_W(Z)}{\partial W_k\partial W_l}\Vert\leq Cte(1+\Vert Z\Vert^2) $ and \\
$\Vert \frac{\partial^2 F_W(Z)}{\partial W_k\partial W_l}-\frac{\partial^2 F_W^0(Z)}{\partial W_k\partial W_l}\Vert\leq Cte\Vert W-W^0 \Vert(1+\Vert Z\Vert^3)$,  by standard arguments found, for example, in \cite{Yao} we get 
\[
\sqrt{n}\left(\hat W_n-W^0\right)\stackrel{Law}{\rightarrow}{\cal N}(0,I^{-1}_0)
\]
\(
\blacksquare
\)

\subsection{Asymptotic distribution of $T_n$}
In this section, we write  $\hat W_n=\arg\min_{W\in \Theta_s}U_n(W)$ and \\
$\hat W^0_n=\arg\min_{W\in \Theta_q} U_n(W)$, where $\Theta_q$ is view as a subset of $\mathbb R^s$. The asymptotic distribution of $T_n$ is then a consequence of the previous section, namely, if we have to replace $nU_n(W)$ by its Taylor expansion around $\hat W_n$ and $\hat W^0_n$, following  \cite{Vandervaart} chapter 16 we have :
\[
T_n=\sqrt{n}\left(\hat W_n-\hat W^0_n\right)^T I_0\sqrt{n}\left(\hat W_n-\hat W^0_n\right)+o_P(1)\stackrel{\cal D}{\rightarrow}\chi^2_{s-q}
\]
\section{Conclusion}
It has been show that, in the case of multidimensional output, the cost function $U_n(W)$ leads to a test for the number of parameters in MLP simpler than with the traditional mean square cost function. In fact the estimator $\hat W_n$ is also more efficient than the least square estimator (see \cite{Rynkiewicz}). We can also remark that $U_n(W)$ matches with twice the ``concentrated Gaussian log-likelihood'' but we have to emphasize, that its nice asymptotic properties need only moment condition on $\varepsilon$ and $Z$, so it works even if the distribution of the noise is not Gaussian. An other solution could be to use an approximation of the covariance error matrix to compute generalized least square estimator :
\[
\frac{1}{n}\sum _{t=1}^{n}\left(Y_{t}-F_{W}\left(Z_{t}\right)\right)^{T}\Gamma ^{-1}\left(Y_{t}-F_{W}\left(Z_{t}\right)\right),
\]
assuming that $\Gamma $ is a good approximation of the true covariance
matrix of the noise $\Gamma(W^0)$.
However it take time to compute a good the matrix $\Gamma$ and if we try to compute the best matrix $\Gamma$ with the data, it leads to the cost function $U_n(W)$ (see for example \cite{Gallant}). 

Finally, as we see in this paper, the computation of the derivatives of $U_n(W)$ is easy, so we can use the effective differential optimization techniques to estimate $\hat W_n$ and numerical examples can be found in \cite{Rynkiewicz}.

\bibliographystyle{asmda2005References}

\bibliography{bibliography-asmda2005}

\begin{thebibliography}{}

\bibitem[\protect\citeauthoryear{Fukumizu}{2003}]{Fukumizu}
K.~Fukumizu.
\newblock Likelihood ratio of unidentifiable models and multilayer neural
  networks.
\newblock {\em Annals of Statistics}, 31:3:533--851, 2003.

\bibitem[\protect\citeauthoryear{Gallant}{1987}]{Gallant}
R.A. Gallant.
\newblock {\em Non linear statistical models}.
\newblock J. Wiley and Sons, New-York, 1987.

\bibitem[\protect\citeauthoryear{Hall and Heyde}{1980}]{Hall}
P.~Hall and C.~Heyde.
\newblock {\em Martingale limit theory and its applications}.
\newblock Academic Press, New-York, 1980.

\bibitem[\protect\citeauthoryear{Magnus and Neudecker}{1988}]{Magnus}
Jan~R. Magnus and Heinz Neudecker.
\newblock {\em Matrix differential calculus with applications in statistics and
  econometrics}.
\newblock J. Wiley and Sons, New-York, 1988.

\bibitem[\protect\citeauthoryear{Rynkiewicz}{2003}]{Rynkiewicz}
J.~Rynkiewicz.
\newblock Estimation of multidimensional regression model with multilayer
  perceptrons.
\newblock In J.~Mira and J.R. Alvarez, editors, {\em Computational methods in
  neural modeling}, volume 2686 of {\em Lectures notes in computer science},
  pages 310--317, 2003.

\bibitem[\protect\citeauthoryear{Sussman}{1992}]{Sussman}
H.J. Sussman.
\newblock Uniqueness of the weights for minimal feedforward nets with a given
  input-output map.
\newblock {\em Neural Networks}, pages 589--593, 1992.

\bibitem[\protect\citeauthoryear{Van~der Vaart}{1998}]{Vandervaart}
A.W. Van~der Vaart.
\newblock {\em Asymptotic statistics}.
\newblock Cambridge University Press, Cambridge, UK, 1998.

\bibitem[\protect\citeauthoryear{Yao}{2000}]{Yao}
J.~Yao.
\newblock On least square estimation for stable nonlinear ar processes.
\newblock {\em The Annals of Institut of Mathematical Statistics}, 52:316--331,
  2000.

\end{thebibliography}

\end{document}